\documentclass[11pt]{amsart}
\usepackage{amsthm,amsmath,amssymb,amscd}

{\end{list}}
{
   \newtheorem{theorem}{Theorem}[section]
   \newtheorem{proposition}[theorem]{Proposition}
   \newtheorem{lemma}[theorem]{Lemma}

   \newtheorem{corollary}[theorem]{Corollary}

}
{\theoremstyle{definition}

   \newtheorem{example}[theorem]{Example}
   \newtheorem{definition}[theorem]{Definition}
}
{\theoremstyle{remark}
   \newtheorem{remark}[theorem]{Remark}

}
\newcommand{\RR}{{\mathbb{R}}}

\newcommand{\NN}{{\mathbb{N}}}

\newcommand{\ZZ}{{\mathbb{Z}}}

\newcommand{\cA}{{\mathcal A}}
\newcommand{\cB}{{\mathcal B}}
\newcommand{\cC}{{\mathcal C}}
\newcommand{\cD}{{\mathcal D}}
\newcommand{\cF}{{\mathcal F}}

\newcommand{\cL}{{\mathcal L}}

\newcommand{\cM}{{\mathcal M}}

\newcommand{\cP}{{\mathcal P}}

\newcommand{\Coker}{\operatorname{Coker}}
\newcommand{\Ker}{\operatorname{Ker}}
\newcommand{\res}{\operatorname{res}}

\newcommand{\Star}{\operatorname{Star}}

\newcommand{\Hom}{{\operatorname{Hom}}}

\newcommand{\isom}{\simeq}

\newcommand{\cbl}{\cite{BL} }
\newcommand{\cbb}{\cite{BBFK} }

\newcommand{\tensor}{\otimes}

\newcommand{\Cd}{C^\bullet}
\newcommand{\Po}{Poincar\'e }
\newcommand{\DD}{{\mathbb{D}}}
\newenvironment{Sb}{_\bgroup%
   \everymath{\scriptstyle}\array{@{}c@{}}%
 }{
   \endarray\egroup
 }
\newcommand{\setmin}{{\smallsetminus}}

\begin{document}
\title{The cd-index of fans and lattices.}

\author{Kalle Karu}
\thanks{The author was partially supported by NSERC}
\address{Mathematics Department\\ University of British Columbia \\
  1984 Mathematics Road\\
Vancouver, B.C. Canada V6T 1Z2}
\email{karu@math.ubc.ca}

\maketitle


\section{Introduction}

The number of flags in a complete fan, or more generally in an
Eulerian poset, is encoded in the cd-index \cite{BB,BK,S}. The goal
of this article 
is to prove the non-negativity of the cd-index for complete fans,
regular CW-spheres \cite{BK} and Gorenstein* posets \cite{S}.

Let us start by recalling the definition of the cd-index for complete
fans. If $\Delta$ is a complete simplicial fan in $\RR^n$, define
$\cA(\Delta)$ to 
be the space of continuous conewise polynomial functions on
$\Delta$ (i.e., functions on $\RR^n$ that restrict to polynomials on each cone
of $\Delta$). Then $\cA(\Delta)$ is a graded free module over the ring $A$ of
global polynomial functions. We write its \Po polynomial 
\[ p_n(\Delta) = \sum_k h_k t^k,\]
where $h_k = \dim (\cA(\Delta) \tensor_A \RR)^k$. When the fan
$\Delta$ is rational, it corresponds to a  toric variety $X(\Delta)$
and the numbers $h_k$ are the even Betti numbers of
the cohomology of $X(\Delta)$. In any case, they satisfy the \Po
duality $h_k = h_{n-k}$. 

When the fan $\Delta$ is complete but not necessarily simplicial, one
can apply the same construction to the first barycentric subdivision
$B(\Delta)$ of $\Delta$. We can label a $1$-dimensional cone of
$B(\Delta)$ according to the dimension of the cone in $\Delta$ that it
is a barycenter of. This gives a multi-grading on $\cA(B(\Delta))$ by $\NN^n$.  
It is also possible to adjust the $A$-module structure of
$\cA(B(\Delta))$ so that it preserves this multigrading. 
The corresponding \Po polynomial is
\[ P_n(\Delta) = \sum_S h_S t^S,\]
where $S\in \NN^n$, $t^S = t_1^{S_1}\cdots t_n^{S_n}$ and $h_S$ is the
dimension of the degree $S$ part of $\cA(B(\Delta)) \tensor_A \RR$. From
the \Po duality $h_S = h_{\overline{S}}$, where $\overline{S} =
(1,\ldots,1)-S$, we see that $h_S$ is zero unless $0\leq S_i \leq
1$. Thus, the sum in the \Po polynomial runs over subsets
$S\subset \{1,\ldots,n\}$. 

It turns out that the \Po polynomial $P_n(\Delta)$ can be written as a
polynomial in non-commuting variables $c$ and $d$ as follows. To a
monomial in $c$ and $d$ we associate a polynomial in $t_1,\ldots,t_n$
by replacing $c$ with $t_i+1$ and $d$ with $t_i+t_{i+1}$, e.g.,
\[ ccdcd = (t_1+1)(t_2+1)(t_3+t_4)(t_5+1)(t_6+t_7).\]
The replacement is done by starting from the left and using each $t_i$
exactly once in increasing order. Define the degree of $c$ to be $1$
and the degree of $d$ to be $2$. Then Bayer and Klapper \cite{BK} show
that $P_n(\Delta)$ is a homogeneous $cd$-polynomial of degree $n$,
called the $cd$-index of $\Delta$.

\begin{example} For any complete $n$-dimensional fan the coefficient
  of $c^n$ in the 
  $cd$-index is $1$ because it is the dimension of the degree $0$ part
  in $\cA(B(\Delta))$. Thus, $c$ is the $cd$-index of the complete
  $1$-dimensional fan. The $cd$-index of the complete $2$-dimensional
  fan with $k$  
  maximal cones is $c^2+(k-2) d$. If $\Delta$ is the $3$-dimensional
  fan over the faces of a pyramid (with square base) then its
  $cd$-index is $c^3+3cd+3dc$.

A simple way to construct examples in low dimension is to find a
subfan $\Pi$ of the $(n-1)$-skeleton $\Delta^{\leq n-1}$ of $\Delta$
such that $\Pi$ contains the $(n-2)$-skeleton  $\Delta^{\leq n-2}$ and
is combinatorially equivalent to a complete $(n-1)$-dimensional
fan. Let $\sigma_1,\ldots,\sigma_m$ be the $(n-1)$-dimensional cones
of $\Delta\setmin \Pi$. Then 
\[ P_n(\Delta) = P_{n-1}(\Pi) c + \sum_i P_{n-2}(\partial \sigma_i) d.\]
If $\Delta$ is the fan over the faces of a $3$-dimensional polytope
then finding such a $\Pi$ amounts to finding a Hamiltonian cycle in
the edge graph of the polytope.

\end{example}

We prove another combinatorial construction of the $cd$-index. Let
$\Delta^{\leq m}$ be the $m$-skeleton of $\Delta$ and 
$I(\Delta^{\leq m})$ the set of functions $f:\Delta^{\leq m}\to \ZZ$.
We define the following operations on such functions. First let $E: I(
\Delta^{\leq m}) \to I( \Delta^{\leq m})$ 
be
\[ E(f)(\sigma) = \sum_{\tau \in \Delta^{\leq m},\; \sigma\leq \tau}
(-1)^{m-\dim\tau} f(\tau).\]
Also, let $C: I( \Delta^{\leq m}) \to I( \Delta^{\leq m-1})$ be the
restriction, and let $D: I( \Delta^{\leq m}) \to I( \Delta^{\leq m-2})$ be
\[ D = C\circ (E-Id) \circ C.\]
Now if $w(c,d)$ is a $cd$-monomial of degree $n$, then applying
$w(C,D)$ to an element $f\in I(\Delta)$ gives and element in
$I(\Delta^0)$, that is, an integer number. 

\begin{proposition}\label{prop-constr} Let $1_\Delta \in I(\Delta)$ be
  the constant function 
  $1$. If $w(c,d)$ is any $cd$-monomial of degree $n$ then
  $w(C,D)(1_\Delta)$ is the coefficient of the monomial $w(c,d)$ in
  the $cd$-index of $\Delta$.
\end{proposition}

To prove that the $cd$-index of a complete fan is non-negative, we
consider sheaves on $\Delta$. For a sheaf of finite
dimensional vector spaces, taking the dimension of stalks gives an
element in $I(\Delta)$ with non-negative integer values. We show
that the operations $C$ and $D$ have their counterparts in the
category of sheaves, hence applying any $cd$-monomial to $1_\Delta$
gives a function with non-negative values.

The $cd$-index can be defined more generally for Eulerian posets.
Bayer and Klapper \cite{BK} conjectured that the $cd$-index of a
regular $CW$-sphere has non-negative coefficients. Stanley \cite{S}
generalized this conjecture to Gorenstein* posets which are defined as
follows. If $\Lambda$ is a rank $n$ graded poset, let $B(\Lambda)$ be
the simplicial complex of chains in $\Lambda\setmin\{{\bf
  0,1}\}$. Then $\Lambda$ is 
called Gorenstein* if $B(\Lambda)$ is a homology manifold with
the homology of the sphere $S^{n-1}$ (see Section~\ref{sec-Gorenstein}
for more on Gorenstein* posets). Complete fans and regular
$CW$-spheres are examples of Gorenstein* posets; in these cases
$|B(\Lambda)|$ is homeomorphic to $S^{n-1}$.

The main result we prove is

\begin{theorem}\label{thm-main} Let $\Lambda$ be a complete fan 
  or a Gorenstein* poset. Then the $cd$-index of $\Lambda$ has
  non-negative integer coefficients.
\end{theorem}

Stanley \cite{S} has shown that there are no other linear inequalities
satisfied by the $cd$-indices of all Gorenstein* posets than the
non-negativity of its coefficients.

Note that it suffices to prove the theorem for Gorenstein* posets
because this case also covers complete fans. Nevertheless we prefer to
use the terminology of fans and prove the theorem for complete fans,
after which we explain which parts of
the proof also apply to Gorenstein* posets.

Since the $cd$-index of a complete fan has non-negative integer
coefficients, one can ask what do these coefficients count? We do not
give an answer to this question. The closest we come to answering it is the
construction using the subfan $\Pi$ above. Suppose in an ideal
world every complete $n$-dimensional fan $\Delta$ has a subfan $\Pi$ of the
$n-1$-skeleton as above. One can think of $\Pi$ as a maximal submanifold
of $\Delta^{\leq n-1}$. If the $cd$-index of $\Delta$ is $fc+gd$ then
$f$ accounts for the ``manifold part'' of  $\Delta^{\leq n-1}$ and $g$
for the rest. Thus, the $cd$-index is a universal invariant that
describes the extraction of a maximal submanifold of $\Delta^{\leq
  n-1}$. This description is rather vague, but the principle of
extracting the manifold part of $\Delta^{\leq n-1}$ also holds when we
describe in Section~\ref{sec-shelling} how to cut $\Delta^{\leq n-1}$
into manifolds with boundary and compute the $cd$-index from these
pieces.

{\bf Acknowledgments.} I would like to thank Louis Billera for
introducing me to the $cd$-index.

\section{Preliminaries}

We will use the terminology of fans \cite{F}. If $\Delta$ is a fan and
$\sigma\in\Delta$ a cone, let $[\sigma]$ be the fan consisting of
all faces of $\sigma$, $\partial\sigma = [\sigma]\setmin\{\sigma\}$,
and $\Star\sigma = \{\gamma\in\Delta| \gamma\geq \sigma\}$. For a graded set the
superscript refers to the degree. If $\Delta$ is a fan or a graded poset, then
$\Delta^{\leq m}$ is the $m$-skeleton of $\Delta$. The dimension of a
cone, or the degree of an element in a poset, is denoted
$\dim(\sigma)$. 

\subsection {Fan spaces}

Let us recall the definition of a topology on a fan and the sheaves on
it \cite{BBFK, BL}. 

Given a fan $\Delta$ (considered as a finite partially ordered set),
the topology on $\Delta$ is generated by subsets $[\sigma]$ for
$\sigma\in\Delta$. In other words, open sets in the topology are the
subfans of $\Delta$. A sheaf $F$ of $\RR$-vector spaces is given by
specifying for 
each $\sigma\in\Delta$ the stalk $F_\sigma$ and for $\tau<\sigma$ the
restriction map $res^\sigma_\tau: F_\sigma\to F_\tau$, satisfying the
compatibility condition $res^\sigma_\tau \circ res^\gamma_\sigma =
res^\gamma_\tau$ for $\gamma>\sigma>\tau$. More abstractly, a sheaf is
a functor from the category of the poset $\Delta$ to the category of
vector spaces. For $U\subset \Delta$ an open set, we denote by $F(U)$
the space of sections of $F$ on $U$. Thus, $F(\Delta) = H^0(\Delta,F)$.

The structure sheaf $\cA$ on $\Delta$ is defined as follows. For
$\sigma\in\Delta^k$, let 
\[ \cA_\sigma = A_k := \RR[x_1,\ldots,x_k],\]
and for $\tau<\sigma$, $\tau\in\Delta^l$, let the restriction map
$res^\sigma_\tau$ be the standard surjection $A_k\to A_l$ (sending
$x_i\mapsto x_i$ for $i\leq l$ and $x_j\mapsto 0$ for $j>l$). Note
that if $\Delta$ has 
dimension $n$ then the sheaf $\cA$ is multi-graded by $\NN^n$.

Let $\cM(\cA)$ be the category of (multi-graded) locally free flabby
$\cA$-modules: a sheaf $\cF$ lies in  $\cM(\cA)$ if the stalks
$\cF_\sigma$ are finitely generated graded free $\cA_\sigma$ modules
and the restrictions 
\[ \cF_\sigma \to \cF(\partial\sigma) \]
are surjective morphisms, compatible with the module structure and
grading. There exists an indecomposable sheaf $\cL \in \cM(\cA)$
satisfying the following conditions:
\begin{enumerate}
\item $\cL_0 = \RR$.
\item The map induced by restriction
\[ \cL_\sigma \tensor_{A_k} \RR \to \cL(\partial\sigma) \tensor_{A_k}
\RR \]
is an isomorphism for all $\sigma\in\Delta^k$. 
\end{enumerate}
The sheaf $\cL$ is unique up to an isomorphism and it can be
constructed by induction on the dimension of the cones using the two
properties. In fact, $\cL(\partial\sigma)$ is a free $A_{k-1}$-module
for any $\sigma\in\Delta^k$ and we can put
\[ \cL_\sigma = \cL(\partial\sigma) \tensor_{A_{k-1}} A_k.\]

Now let $B(\Delta)$ be a barycentric subdivision of $\Delta$. It is
constructed by choosing for each $\sigma\in\Delta$ a vector $v_\sigma$
in the relative interior of $\sigma$ and star-subdividing the fan at
$v_\sigma$, starting with cones of maximal dimension, then one
dimension smaller, and so on. The fan $B(\Delta)$ is simplicial and its
cones are in one-to-one correspondence with chains of cones $0<
\sigma_1<\sigma_2<\ldots<\sigma_k$ in $\Delta$. This gives $B(\Delta)$ a
multi-grading by subsets $S = \{\dim\sigma_1,\ldots,\dim\sigma_k\}
  \subset \{1,\ldots,n\}$. 

Let us define a sheaf of rings $\cB$ on $B(\Delta)$. If $\sigma\in
B(\Delta)^S$, set
\[ \cB_\sigma = B_S := \RR[x_i]_{i\in S},\]
and let the restriction maps $res^\sigma_\tau$ be the standard
surjections $B_S \to B_T$. This sheaf again
is multi-graded by $\NN^n$. We can define the indecomposable
$\cB$-module $\cL$ as before, but since $\cB$ itself is flabby, it
follows that $\cL \isom \cB$. To avoid confusion, we will use the
letter $\cL$ for the indecomposable $\cA$-module on $\Delta$ and the
letter $\cB$ for the similar $\cB$-module on $B(\Delta)$.  

It should be noted that the definition of topology, structure sheaf,
and the sheaf $\cL$ only depend on the partially ordered set $\Delta$,
and not its realization as a fan. We say that two fans are
combinatorially equivalent if their posets are
isomorphic. Combinatorially  equivalent fans give rise to isomorphic
theories of sheaves on them. 

Define a map of posets $\pi: B(\Delta)\to\Delta$ by sending a cone in
$B(\Delta)$ to the smallest cone in $\Delta$ containing it.
If $\sigma\in  B(\Delta)$ corresponds to a chain $0<
\sigma_1<\sigma_2<\ldots < \sigma_k$, then $\pi(\sigma) = \sigma_k$.
 This map is continuous in the fan topologies. We can also extend it
 to a map of ringed spaces, letting
\[ \cA_{\pi(\sigma)} = \cA_{\sigma_k} = A_{\dim\sigma_k} \to
  B_{\{\dim\sigma_1,\ldots,\dim\sigma_k\}} = \cB_\sigma \]
be the standard projection. 

\begin{lemma}\label{lem-iso} There exists an isomorphism of sheaves of
  $\cA$-modules on $\Delta$ 
\[ \cL \isom \pi_* \cB.\]
\end{lemma}

{\bf Proof.} Let us check that the sheaf $\pi_* \cB$ satisfies the two
conditions in the 
definition of $\cL$. The first condition is trivially true. For the
second note that combinatorially for any $\sigma\in \Delta^k$
\[ B([\sigma]) \isom B(\partial\sigma) \times [\rho],\]
where $\rho$ is a $1$-dimensional cone. It follows that 
\[ \pi_* \cB ([\sigma]) = \cB(B([\sigma])) \isom \cB(B(\partial\sigma))
\otimes_{A_{k-1}} A_k.\]
As we will see below, $\cB(B(\partial\sigma))$ is a free
$A_{k-1}$-module, so the sheaf $\pi_* \cB$ lies in
$\cM(\cA)$, hence is isomorphic to $\cL$. \qed

\subsection{Duality}

Let us summarize some of the duality results proved in
\cite{BBFK,BL}. The structure sheaves defined in those papers 
differ from the ones used here, but the constructions are general
enough to work also for sheaves $\cA$ and $\cB$.
We will mostly be using the category-theoretic constructions from \cite{BL}.  
 
Fix an $n$-dimensional complete fan $\Delta$.

\begin{enumerate}
\item Let $F$ be a sheaf on
 $\Delta$. The cellular complex of $F$ is
\[ \Cd_n(\Delta, F) =  0\to C^0 \to C^1\to \cdots \to C^n\to
0,\]
where 
\[ C^i = \bigoplus_{\sigma\in \Delta^{n-i}} F_\sigma,\]
and the differentials $\delta: C^i\to C^{i+1}$ are defined by choosing an
orientation for each cone $\sigma\in\Delta$ and setting
$F_\sigma\to F_\tau$ to be $\pm res^\sigma_\tau$ depending on whether the
orientations agree or not. 

On a complete fan $\Delta$ the sheaf cohomology of $F$ can be computed
using the cellular complex: 
\[ H^i(\Delta,F) \isom H^i(\Cd_n(\Delta, F)).\]

\item Since the sheaf $\cL$ is flabby, its nonzero cohomology
  vanishes. Moreover, $H^0(\Delta,\cL)$ is a graded free $A_n$-module
  of finite rank. Similar result holds for the sheaf $\cB$ on $B(\Delta)$, and
  from Lemma~\ref{lem-iso} we get 
\[ H^0(\Delta,\cL) \isom  H^0(B(\Delta),\cB).\]
Denote the \Po polynomial 
\[ P_n(\Delta) = \sum_{S\in \NN^n} h_S t^S, \]
where $h_S = \dim (H^0(\Delta,\cL)\tensor_{A_n} \RR)^S$ and $t^S =
t_1^{S_1} \cdots t_n^{S_n}$. 

\item Let $\omega_n = A_n[(-1,\ldots,-1)]$. I.e., $\omega_n$ is the
free $A_n$-module of rank $1$ placed in degree $(1,1,\ldots,1)$. Then 
\[ Hom_{A_n}(H^0(\Delta,\cL),\omega_n) \isom H^0(\Delta,\cL).\]
In particular, the numbers $h_S$ satisfy \Po duality
\[ h_S = h_{\overline{S}},\]
where $\overline{S} = (1,1,\ldots,1)-S$. It follows that we may index
the numbers 
$h_S$ by subsets $S\subset\{1,\ldots,n\}$.

\item More generally, consider the $m$-skeleton $\Delta^{\leq m}$ of
  $\Delta$. The sheaf cohomology of $\cL$ on $\Delta^{\leq m}$
  vanishes in nonzero degrees because $\cL$ is flabby, and it is a free
  $A_m$-module in degree $0$. There is an  isomorphism in the derived
  category $D^b(\cA)$ of finitely generated $\cA$-modules
\[ R\Gamma( \Delta^{\leq m}, \cL) \isom RHom (\Cd_m(\Delta^{\leq m},
  \cL),\omega_m).\]
In particular, the cellular complex of $\cL$ on $\Delta^{\leq m}$ has
no higher cohomology; the degree zero cohomology is a free
$A_m$-module, dual to the module of global sections. When $m=n$ this
together with $(1)$ gives \Po duality.

\end{enumerate}

The cellular complex can be used to express the numbers $h_S$ in terms
of flag numbers of $\Delta$. First, let us generalize the definition
of \Po polynomial. 

If $M$ is a graded free $A_n$ module, we define its \Po polynomial
$P_n(M)$ as above. Since the module may have elements in negative
degree, $P_n(M)$ is a Laurent polynomial in general.  For a complex of
finite length $M^\bullet$ of free $A_n$-modules we 
define the \Po polynomial by taking the alternating sum 
\[ P_n(M^\bullet) = \sum_i (-1)^i P_n(M^i).\]
Finally, for any finitely generated $A_n$-module (or a bounded complex
of such modules) we first choose a free 
resolution and then compute its \Po polynomial. It is easy to see that
this is well defined and in fact it gives the \Po polynomial of an
object in $D^b(A_n)$. Notice also that the \Po polynomial of the dual
module $Hom_{A_n}(M,A_n)$ is obtained from the \Po polynomial of $M$
by replacing $t_i$ with $t_i^{-1}$.

Since the cellular complex $\Cd_n(B(\Delta),\cB)$ is a resolution of
$H^0(B(\Delta),\cB) \isom H^0(\Delta,\cL)$, we have
\[ P_n(\Delta) = \sum_{\sigma\in B(\Delta)} (-1)^{n-\dim\sigma} P_n
(\cB_\sigma).\]
The $A_n$-module $\cB_\sigma = B_S$ has the Koszul resolution
\[ \ldots \to \Lambda^2 V \otimes A_n \to V\otimes A_n \to A_n \to
B_S,\]
where $V$ is a graded vector space with basis $\{x_i\}_{i\notin S}$. It
follows that if $\sigma$ has multi-degree $S$ then 
\[ P_n(\cB_\sigma) = \prod_{i\notin S} (1-t_i).\]
Thus
\[  P_n(\Delta) = \sum_S \sum_{\sigma\in B(\Delta)^S} (-1)^{n-|S|}  
 \prod_{i\notin S} (1-t_i).\]
The numbers $f_S = |B(\Delta)^S|$ are called flag numbers of $\Delta$
-- they count the number of chains (or flags) in $\Delta$ of type
$S$. We can express the coefficient of $t^T$ in the formula above as
\[ h_T =  \sum_{S\subset \overline{T}} (-1)^{n-|S|-|T|} f_S.\]
Since $h_T = h_{\overline{T}}$, the sum can also be written as
\[ h_T =   \sum_{S\subset T} (-1)^{|T\setmin S|} f_S.\]
This last formula is often used to define the numbers $h_T$. Note that
we can also invert this formula to express $f_S$ in terms of $h_T$.

\subsection{Existence of the cd-index}

Let $c$ and $d$ be non-commuting variables of degree $1$ and $2$,
respectively, and let $\RR[c,d]$ be the polynomial ring in $c$ and
$d$. We define an injective $\RR$-vector space homomorphism $\phi:
\RR[c,d] \to \RR[t_1, t_2,\ldots]$  
by induction on degree. If $f(c,d)$ is a homogeneous polynomial of
degree $k$ on which $\phi$ is defined, let
\begin{align*} 
\phi(f(c,d) c) & = \phi(f(c,d)) (t_{k+1}+1),\\
 \phi(f(c,d) d) & = \phi(f(c,d)) (t_{k+1}+t_{k+2}).
\end{align*}
We say that $g(t_1,\ldots,t_n)$ is a $cd$-polynomial if it
lies in the image of $\phi$. When talking about the degree of a
$cd$-polynomial we always mean the degree in $c$ and $d$.

The following proposition was proved by Fine, Bayer and Klapper
\cite{BK}. The proof we give is essentially the same as in 
\cite{BK, S}.

\begin{proposition}\label{cd-exist}
Let $\Delta$ be a complete $n$-dimensional fan. Then $P_n(\Delta)$ is
a homogeneous $cd$-polynomial of degree $n$.
\end{proposition}

{\bf Proof.} By \Po duality
\[ 2P_n(\Delta) = P_n(H^0(\Delta,\cL)) +
P_n(Hom_{A_n}(H^0(\Delta,\cL),\omega_n)).\]
Since the cellular complex of $\cL$ is a resolution of
$H^0(\Delta,\cL)$, we have
\begin{equation}  2P_n(\Delta) = \sum_{\sigma\in\Delta}
  (-1)^{n-\dim\sigma} ( P_n (\cL_\sigma) + P_n
  (RHom_{D^b(A_n)}(\cL_\sigma,\omega_n))), \label{sum1}
\end{equation}
where $RHom_{D^b(A_n)}(\cL_\sigma,\omega_n)$ is the derived $Hom_{A_n}$,
computed as the $Hom_{A_n}$ of a free resolution of $\cL_\sigma$.
It suffices to prove that each term in the sum above is a $cd$-polynomial.

Let us fix a $\sigma\in\Delta^k$, $k>0$. By induction on dimension we may
assume that $P_k(\cL_\sigma) = P_{k-1}(\partial\sigma)$ is a homogeneous
$cd$-polynomial of degree $k-1$. From the Koszul resolution we get
\[  P_n (\cL_\sigma) = P_k (\cL_\sigma) \prod_{i>k}(1-t_i).\]
\Po duality applied to $\partial\sigma$ implies that
\[P_k (Hom_{A_k}(\cL_\sigma,\omega_k)) = 
P_k (\cL_\sigma)(t_1^{-1},\ldots,t_k^{-1}) t_1\cdots t_k  =
P_k (\cL_\sigma) t_k.\]
From this we compute
\begin{align*} P_n (RHom_{D^b(A_n)}(\cL_\sigma,\omega_n)) 
&= P_n (\cL_\sigma)(t_1^{-1},\ldots,t_n^{-1}) t_1\cdots t_n \\
&= P_k (\cL_\sigma)(t_1^{-1},\ldots,t_k^{-1}) \prod_{i>k}(1-t_i^{-1})
  t_1\cdots t_n \\ 
&= P_k (\cL_\sigma) t_k \prod_{i>k}(t_i-1) \\
& = (-1)^{n-k}P_k (\cL_\sigma) t_k  \prod_{i>k}(1-t_i).
\end{align*}
Thus, a term corresponding to $\sigma$ in the sum~(\ref{sum1}) is
\[ (-1)^{n-k} P_k (\cL_\sigma) ( 1+(-1)^{n-k} t_k) \prod_{i>k}(1-t_i)
= P_k (\cL_\sigma) ( (-1)^{n-k} + t_k) \prod_{i>k}(1-t_i).\]
Since $(1-t_i)(1-t_{i+1}) =(1+t_i)(1+t_{i+1})-2(t_i+ t_{i+1}) =
c^2-2d$, we get
\[ P_k (\cL_\sigma) ( (-1)^{n-k} +  t_k) \prod_{i>k}(1-t_i) =
\begin{cases} P_k (\cL_\sigma) c (c^2-2d)^{\frac{n-k}{2}} \quad
  \text{if $n-k$ is even} \\
-P_k (\cL_\sigma) (c^2-2d)^{\frac{n-k+1}{2}} \quad
  \text{if $n-k$ is odd}.
\end{cases} \]
When $\sigma$ is the zero cone then a similar computation shows that
the contribution from $\sigma$ to the sum~(\ref{sum1}) is
\begin{equation} \label{zero} ((-1)^n+1) \prod_{i>1}^n(1-t_i) = \begin{cases}
  2(c^2-2d)^{\frac{n}{2}} & \text{ if $n$ is even} \\
0 & \text{ if $n$ is odd}.
\end{cases} 
\end{equation}
\qed 

\begin{remark} 
\begin{enumerate}
\item The proof above  gives a formula of Stanley \cite{S} for recursively
computing the $cd$-index:
\[ P_n(\Delta) = \frac{1}{2} [ \sum\begin{Sb} \sigma\in\Delta^k\setmin\{0\}\\ n-k
    \text{ even} \end{Sb} P_{k-1} (\partial\sigma) c (c^2-2d)^{\frac{n-k}{2}} -  
 \sum \begin{Sb}\sigma\in\Delta^k\setmin\{0\}\\ n-k
    \text{ odd}\end{Sb} P_{k-1} (\partial\sigma) (c^2-2d)^{\frac{n-k+1}{2}}
 +\epsilon],\]
where $\epsilon$ is the $cd$-polynomial (\ref{zero}) above.
\item  The proof also shows that for any fan
 $\Delta$, not necessarily complete, the sum of the \Po polynomials of
 the cellular complex of $\cL$ and its dual is a $cd$-polynomial. 
\end{enumerate}
\end{remark}

\subsection{Gorenstein* posets}
\label{sec-Gorenstein}

It is clear that the topology and the ringed space structure on a fan
$\Delta$ only depend on the poset of $\Delta$. Let us
discuss the conditions a poset must satisfy in order to apply the
results of this section. Gorenstein* posets defined below  form the
largest class of posets for which everything works as for fans.

We refer to \cite{S1} for the terminology and results on posets.
By a graded rank $n$ poset we mean a finite graded poset $\Lambda$
with a unique 
minimal element $\bf{0}$ in degree $0$ and a unique maximal element
$\bf{1}$ in degree $n+1$. If such elements do not exist, we add
them. The poset $\Lambda$ being graded in particular implies that all
maximal chains have the same length.

For any graded poset the definition of topology, structure sheaf $\cA$ and
the irreducible sheaf $\cL$ goes through as for fans. Note that by a
sheaf on a poset $\Lambda$ we mean a sheaf on  $\Lambda \setmin
\{{\bf 1}\}$. What we don't have for an arbitrary poset $\Lambda$ is
the cellular complex of a sheaf on $\Lambda$. To define the cellular
complex, we need both the restriction maps $\res^\sigma_\tau$ and the
orientations $or^\sigma_\tau = \pm 1$.

Let $\Lambda$ be a rank $n$ graded poset. Define $\Lambda$ to be 
{\em orientable} if for every $x<y<{\bf 1}$, where $y$ covers $x$, we
can choose an orientation  
\[ or^y_x = \pm 1,\]
such that the sequence $\Cd_n(\Lambda,\RR_\Lambda)$ for the constant
sheaf $\RR_\Lambda$ defined the same way as for fans above is a
complex, $\delta\circ \delta= 0$. Note that we do not specify
$or^y_x$ for $y=\bf{1}$ because the maximal element $\bf{1}$ does not
appear in the cellular complex.

We say that an orientable poset $\Lambda$ is {\em Gorenstein} if 
\begin{equation}\label{eq-gor} \dim H^i(\Cd_n(\Lambda,\RR_\Lambda)) =
  \begin{cases} 1 &  \text{if 
    $i=0$}, \\ 0 & \text{otherwise} \end{cases}
\end{equation}
and recursively, every proper sub-interval $[x,y] \subset \Lambda$ is
also Gorenstein. (This is not a standard terminology.) Recall that a
graded poset $\Lambda$ is {\em Eulerian} 
if every interval $[x,y]$ for $x<y$ has the same number of elements of
even and odd degree. A Gorenstein poset is clearly Eulerian: the Euler
characteristic of the cellular complex on any half-open interval
$[x,y)$ is $1$.

To apply the duality theory of Bressler and Lunts \cite{BL}, we need a
Gorenstein poset that is actually a lattice.
A complete fan is an example of a Gorenstein lattice. Regular
$CW$-spheres \cite{BK} are Gorenstein posets but they are not usually
lattices. 

The proof of duality in \cite{BL} for fans also involves
the following technicality. If $U$ is a subfan of a complete fan
$\Delta$, define its complement $U'=  \{\gamma\in\Delta: \gamma\cap|U|
= {\bf 0}\}$. Then for the subfan $[\sigma]\subset \Delta$  we have
$[\sigma]'' = [\sigma]$. Also, if $V_\sigma = \Delta\setmin [\sigma]'$
then $V_\sigma$ (or the corresponding space) has a deformation
retraction to $\sigma$. It follows that the cellular complex $\Cd_n$
applied to the constant sheaf $\RR$ on $V_\sigma$ is acyclic in
degrees $<n-1$. Let us call a graded oriented lattice $\Lambda$ {\em
  well-complemented} if it satisfies these two conditions, and
recursively, every proper sub-lattice $[{\bf 0},x]$ is also
well-complemented.   
Now the duality theory holds for well-complemented Gorenstein
lattices. In fact, everything we prove for complete fans also works
word-by-word for such lattices.

Examples of well-complemented lattices are lattices coming from
polyhedral complexes, in particular simplicial complexes such as
$B(\Lambda)$ below. Indeed, for polyhedral complexes the second
condition of well-complemented is satisfied because we have a
deformation retract as in the case of fans. The first condition is
equivalent to the statement that if the vertices of a polyhedron $Q$
are also vertices of $P$ then $Q$ is a face of $P$.

To define the more general notion of a {\em Gorenstein*} poset \cite{S1},
let $\Lambda$ be a graded poset of rank $n$ and let
$B(\Lambda)$ be the lattice of chains in $\Lambda\setmin\{{\bf
0,1}\}$. Then $B(\Lambda)\setmin \{{\bf 0, 1}\}$ is a simplicial
complex with well-defined boundary maps. In particular, $B(\Lambda)$
is orientable. Define $\Lambda$ to be Gorenstein* if
$B(\Lambda)$ is Gorenstein. From the discussion above it
follows that the lattice $B(\Lambda)$ is well-complemented, thus if
$\Lambda$ is Gorenstein*, we can apply the duality theory to $B(\Lambda)$. 

Our definition of a Gorenstein* poset differs from the one
in \cite{S1}, but the two are equivalent. In \cite{S1} a rank $n$
graded poset  
$\Lambda$ is Gorenstein* if (i) the simplicial complex
$B(\Lambda)\setmin \{{\bf 0, 1}\}$  is a homology-manifold, i.e., the
link of every simplex has the the homology of a sphere of appropriate
dimension, and (ii) the entire complex has the homology of the sphere
$S^{n-1}$. These two conditions are equivalent to the cohomology
conditions $(\ref{eq-gor})$ for all intervals $[x,{\bf 1}]$ in
$B(\Lambda)$. For $y\neq {\bf 1}$, the interval $[x,y]$ is the
face lattice of a simplex, hence it satisfies condition
$(\ref{eq-gor})$. It is also not hard to see that if $\Lambda$ is a
Gorenstein* poset then any interval $[x,y]\subset \Lambda$ is again
Gorenstein*. 

The results proved in this section also hold if we replace the fan
$\Delta$ with a Gorenstein* poset $\Lambda$. We define the \Po
polynomial of $\Lambda$ using the sheaf $\cB$ on $B(\Lambda)$ rather
than the sheaf $\cL$ on $\Lambda$. Then  we know that the module of
global sections of $\cB$ satisfies \Po duality.

To prove Proposition~\ref{cd-exist}, we work with the cellular complex
$\Cd_n(B(\Lambda),\cB)$ instead of  $\Cd_n(\Lambda,\cL)$. Write 
\begin{align*}
P_n(\Lambda) = P_n(\Cd_n(B(\Lambda),\cB)) &= \sum_{x\in B(\Lambda)}
(-1)^{n-\dim x} P_n(\cB_x) \\
&= \sum_{\sigma\in\Lambda^k} (-1)^{n-k} \sum_{x\in B(\Lambda), \;
  \pi(x)=\sigma} (-1)^{k-\dim x} P_k(\cB_x) \prod_{i>k}(1-t_i) \\
&= \sum_{\sigma\in\Lambda^k} (-1)^{n-k} P_{k-1}([0,\sigma])
\prod_{i>k}(1-t_i).
\end{align*}
To the term $P_{k-1}([0,\sigma])$ (which is equal to
$P_{k-1}(\partial\sigma)$ in case of fans) we apply induction
assumption and proceed as in the proof of  Proposition~\ref{cd-exist}.

\section{Shelling}\label{sec-shelling}

Stanley \cite{S} has shown that for certain shellable Eulerian posets
the $cd$-index is non-negative. We explain this result for fans. Since
this section is not essential for the general non-negativity proof, we
do not generalize it from fans to posets.

\subsection{Fans with boundary}

Barthel, Brasselet, Fieseler and Kaup studied fans with
boundary in  \cite{BBFK}. A {\em quasi-convex} fan is a
full-dimensional fan $\Delta$ in $\RR^n$ such that the boundary of 
$\Delta$ is a homology manifold. In other words, the lattice of the
boundary fan $\partial\Delta$ is Gorenstein. We also call a fan
quasi-convex if it is  combinatorially equivalent to a quasi-convex
fan as defined above.
In Stanley's terminology a quasi-convex fan is {\em near-Eulerian}. It
can be completed by adding an element with boundary $\partial
\Delta$. Such a completion is a Gorenstein* poset, not necessarily an
actual fan or even a lattice.

It is shown in \cbb that if $\Delta$ is a quasi-convex $n$-dimensional
fan then $\cL(\Delta)$ is a free $A_n$-module. Its dual is the
submodule $\cL(\Delta,\partial\Delta)$ of sections vanishing on the
boundary. The cellular complex of $\cL$ has no non-zero cohomology and 
\[ H^0(\Cd_n(\Delta,\cL)) = \cL(\Delta,\partial\Delta).\]
The module of global sections $\cL(\Delta)$ can also be expressed
using the cellular complex:
\[ \cL(\Delta) = H^0(\Cd_n(\Delta\setmin\partial\Delta,\cL)).\]
We denote by $P_n(\Delta)$ the \Po polynomial of $\cL(\Delta)$ as for
complete fans.

\begin{lemma} Let $\Delta$ be a quasi-convex fan of dimension
  $n$. Then
\[ P_n(\Delta) = f(c,d) + P_{n-1}(\partial\Delta)\]
for some homogeneous $cd$-polynomial $f$ of degree $n$.
\end{lemma}

{\bf Proof.} We complete $\Delta$ to $\overline{\Delta}$ by adding an
element $\sigma$ such that $\partial\sigma=\partial\Delta$. Since
$\cL$ is flabby, we get an exact sequence of free $A_n$-modules
\[ 0\to \cL([\sigma],\partial\sigma) \to \cL(\overline{\Delta}) \to
\cL(\Delta) \to 0.\]
Now $P_n(\overline{\Delta}) = f_1(c,d)$ is a homogeneous $cd$-polynomial
of degree $n$ and $P_n(\cL([\sigma],\partial\sigma)) =
  P_{n-1}(\partial\sigma) t_n$. Hence
\[ P_n(\Delta) = f_1(c,d) -  P_{n-1}(\partial\Delta) t_n = (f_1(c,d) -
P_{n-1}(\partial\Delta)c) + P_{n-1}(\partial\Delta). \qed\]

Let $P_n(\Delta,\partial\Delta) =
P_n(\cL(\Delta,\partial\Delta))$. Then from the lemma and \Po duality
we get
\[  P_n(\Delta,\partial\Delta) =  f(c,d) + P_{n-1}(\partial\Delta)
t_n.\]

Let $\Delta$ be a complete fan and $\Delta^{\leq m}$ its
$m$-skeleton. Then $\cL(\Delta^{\leq m})$ is a free $A_m$-module and
its \Po polynomial is
\[ P_m(\Delta^{\leq m}) = P_n(\Delta)|_{t_{m+1}=\ldots=t_n=0}.\]
It is easy to see that $P_m(\Delta^{\leq m}) = f(c,d)+g(c,d)t_m$ for some
homogeneous $cd$-polynomials $f$ and $g$ of degree $m$ and $m-1$,
respectively. For example, if $m=n-1$ and $P_n(\Delta) =
f(c,d)c+g(c,d)d$, then setting $t_n=0$, we have
\[ P_{n-1}(\Delta^{\leq n-1}) = f(c,d)+g(c,d)t_{n-1}.\]

\begin{lemma}\label{lem-shel} Let $\Delta$ be a quasi-convex
  $n$-dimensional fan and let $\Delta_i$ be quasi-convex 
  $m$-dimensional fans, $m<n$,  
  embedded in the $m$-skeleton  $\Delta^{\leq m}$ of $\Delta$, such
  that every $\sigma\in \Delta^{\leq m}$ occurs in $\Delta_i
  \setmin\partial\Delta_i$ for precisely one $i$. Then 
\[ P_m( \Delta^{\leq m}) = \sum_i P_m(\Delta_i,\partial\Delta_i).\]
\end{lemma}

{\bf Proof.} According to \cite{BL}, the cellular complex $\Cd_m
(\Delta^{\leq m},\cL)$ computes the dual of $\cL(\Delta)$. It
suffices to prove that the \Po polynomials of the dual modules are
equal:
\[ P_m( \Cd_m (\Delta^{\leq m},\cL)) = \sum_i P_m(\Delta_i).\]
This follows from the decomposition of $\Delta^{\leq m}$ into a disjoint union 
\[ \Delta^{\leq m}  = \cup_i (\Delta_i\setmin\partial\Delta_i). \qed\]

\begin{definition} An $n$-dimensional quasi-convex fan $\Delta$ is {\em
    shellable} if there   exists an ordering of the maximal cones
  $\sigma_1,\sigma_2,\ldots,\sigma_N$, such that  $\partial \Delta$
  as well as
\[ \sigma_i^- := [\sigma_i]\cap ([\sigma_1]\cup\ldots\cup
  [\sigma_{i-1}]) \subset \partial\sigma_i\] 
for all  $i=2,\ldots,N$ 
are $(n-1)$-dimensional shellable quasi-convex fans.
\end{definition}

\begin{example}
\begin{enumerate}
\item If $\Delta$ is a complete shellable fan, then $\Delta_i = \sigma_i^-
  \subset  \Delta^{\leq n-1}$ satisfy the conditions of
  Lemma~\ref{lem-shel} (the $(n-1)$-skeleton of $\Delta$ can be
  constructed by starting with $\sigma_N^- = \partial\sigma_N$ and
  attaching to it cells
  $\sigma_{N-1}^-,\sigma_{N-2}^-,\ldots,\sigma_{2}^-$), hence  
\[ P_{n-1}(\Delta^{\leq n-1}) = \sum_{i=2}^N P_{n-1}(\sigma_i^-,\partial
  \sigma_i^-).\]
Let $P_{n-1}(\sigma_i^-,\partial \sigma_i^-) = f_i(c,d) +g_i(c,d)t_{n-1}$
for $cd$-polynomials $f_i$ and $g_i$ of degree $n-1$ and $n-2$,
respectively. Then 
\[ P_n(\Delta) = \sum_{i=2}^N f_i c +g_i d.\]

\item As in the introduction, suppose $\Delta$ is a complete
  $n$-dimensional fan and we have 
\[ \Delta^{n-2} \subset \Pi \subset \Delta^{n-1},\]
where the fan $\Pi$ is combinatorially equivalent to a complete
$(n-1)$-dimensional fan. Let $\Delta^{n-1}\setmin\Pi =
\{\sigma_1,\ldots,\sigma_k\}$. Then $\Pi,
[\sigma_1],\ldots,[\sigma_k]$ provide the $\Delta_i$ in
Lemma~\ref{lem-shel}, hence
\[P_{n-1}(\Delta^{\leq n-1}) = P_{n-1}(\Pi) + \sum_i
P_{n-2}(\partial\sigma_i) t_{n-1}.\]
 From this it follows that
\[ P_n(\Delta) = P_{n-1}(\Pi)c + \sum_i P_{n-2}(\partial\sigma_i)d.\]
\end{enumerate}
\end{example}

\begin{proposition} Let $\Delta$ be a quasi-convex shellable
  $n$-dimensional fan. Then $P_n(\Delta)$ is a $cd$-polynomial with
  non-negative integer coefficients.
\end{proposition} 

{\bf Proof.} Let us fix a shelling $\sigma_1,\ldots,\sigma_N$ of
$\Delta$ and complete $\Delta$ to $\overline{\Delta}$ by adding an
element $\sigma_{N+1}$ such that $\partial\sigma_{N+1} =
\partial\Delta$. (The case when $\partial \Delta = \emptyset$ can be
handled similarly and is left to the reader.) Then
$\sigma_1,\ldots,\sigma_N, \sigma_{N+1}$ is a shelling of
$\overline{\Delta}$ and as in the previous example, we get
\[ P_n(\overline{\Delta}) = \sum_{i=2}^{N+1} f_i c + g_i d,\]
where $P_{n-1}(\sigma_i^-,\partial\sigma_i^-) = f_i(c,d)+g_i(c,d)$; also
$f_{N+1} = P_{n-1}(\partial\Sigma)$ and $g_{N+1} = 0$.
From the exact sequence
\[ 0\to \cL([\sigma_{N+1}],\partial\sigma_{N+1}) \to \cL(\overline{\Delta})
\to \cL(\Delta) \to 0\]
we get
\begin{align*}
 P_n(\Delta) &=   P_n(\overline{\Delta}) -
 P_n([\sigma_{N+1}],\partial\sigma_{N+1}) \\
&= P_n(\overline{\Delta}) - f_{N+1} t_n \\
&= \sum_{i=1}^{N} f_i c + g_i d +P_{n-1}(\partial\Delta).
\end{align*}
By induction on dimension we know that $f_i, g_i$ and
$P_{n-1}(\partial\Delta)$ are homogeneous $cd$-polynomials with
integer coefficients, and so is $P_n(\Delta)$, \qed

\section{Non-negativity of the cd-index}

In this section we prove the main theorem. An essential step in the
proof is to compare $H^0(\Delta^{\leq m},\cL)$ to its dual module
$H^0(\Cd_m(\Delta^{\leq m},\cL))$. To do this we look for a sheaf $F$ on
$\Delta^{\leq m}$ such that the cellular complex of $F\otimes \cL$
computes  $H^0(\Delta^{\leq m},\cL)$. 

Throughout this section we fix a complete $n$-dimensional fan $\Delta$
and an orientation on the cones of $\Delta$ (in fact, we only need the
signs $or^\sigma_\tau$ for $\sigma>\tau$, $\dim\sigma=\dim\tau+1$). We
consider sheaves on the $m$-skeleton $\Delta^{\leq m}$ of $\Delta$.

\begin{definition}\label{def-gor} A sheaf $F$ of finite dimensional vector spaces on
  $\Delta^{\leq m}$ is called {\em semi-Gorenstein} if for any
  $\sigma\in \Delta^{\leq m}$ and any $i>0$
\[ H^i(\Cd_m(\Star \sigma, F)) = 0.\]
If, moreover, for any $\sigma\in \Delta^{\leq m}$
\[ \dim H^0(\Cd_m(\Star \sigma, F)) = \dim F_\sigma,\]
the sheaf $F$ is called {\em Gorenstein}.
\end{definition}

\begin{example}
\begin{enumerate}
\item The constant sheaf $\RR_\Delta$ is Gorenstein on $\Delta$.
\item if $\Delta'\subset \Delta^{\leq m}$ is an $m$-dimensional
  Gorenstein subfan and $\RR_{\Delta'}$ the 
  constant sheaf on $\Delta'$ extended by zero, then $\RR_{\Delta'}$ is 
is a Gorenstein sheaf on $\Delta^{\leq m}$.
\item If $F$ is a semi-Gorenstein sheaf on $\Delta^{\leq m}$ then its
  restriction to $\Delta^{\leq k}$ for $k\leq m$ is again semi-Gorenstein.
\end{enumerate}
\end{example}

\begin{remark} 
From the definition it is clear that a semi-Gorenstein sheaf on
$\Delta^{\leq m}$ is determined by its restriction to the degree
$[m-1,m]$ part $\Delta^{[m-1,m]}$ of $\Delta$. Indeed, if $F$ is
defined on $\Star\sigma \setmin \{\sigma\}$, where $\dim\sigma\leq
m-2$, then 
\[ F_\sigma\isom  H^{m-\dim\sigma-1} (\Cd_m(\Star\sigma
\setmin\{\sigma\},F)),\]
with the obvious boundary maps. However, not every sheaf on
$\Delta^{[m-1,m]}$ extends to a semi-Gorenstein sheaf on $\Delta^{\leq
  m}$.

Let $F,G$ be semi-Gorenstein sheaves on $\Delta^{\leq m}$. A morphism
$F|_{\Delta^{[m-1,m]}} \to 
G|_{\Delta^{[m-1,m]}}$ extends uniquely to a morphism $F\to G$
(inductively, $F_\sigma\to G_\sigma$ is defined as a morphism of
cohomologies above). Since for 
$\tau\in\Delta^{m-1}$ the restriction map 
\[ \oplus_{\sigma>\tau} F_\sigma \to F_\tau \]
is surjective, a morphism $F\to G$ is uniquely determined by its
restriction to $\Delta^m$, but not every morphism  $F|_{\Delta^{m}} \to
G|_{\Delta^{m}}$ extends to a morphism $F\to G$.
\end{remark}

Given a semi-Gorenstein sheaf $F$ on $\Delta^{\leq m}$, we define a new
sheaf $F^\vee$ on $\Delta^{\leq m}$ as follows. Set 
\[ F_\sigma^\vee = H^0(\Cd_m(\Star\sigma,F))^* \]
(the vector space dual), and let the restriction maps be induced by
the projections
\[ \Cd_m(\Star\tau,F) \to \Cd_m(\Star\sigma,F)\]
for $\tau<\sigma$.

A morphism of semi-Gorenstein sheaves $F\to G$ induces a morphism
$G^\vee \to F^\vee$, hence we get a contravariant functor from the
category of semi-Gorenstein sheaves to the category of sheaves.
This functor is exact:  given an exact sequence
of semi-Gorenstein sheaves, applying the functor $(\cdot)^\vee$ gives an
exact sequence (this follows from the vanishing of higher
cohomology: we get an exact sequence of degree zero cohomology).

\begin{lemma} If $F$ is semi-Gorenstein on $\Delta^{\leq m}$, then $F^\vee$ is also
  semi-Gorenstein.
\end{lemma}

{\bf Proof.}
For $\sigma\in\Delta^m$, let $G_\sigma$ be the constant sheaf
$F_\sigma$ on $\{\sigma\}$, extended by zero. Consider the exact
sequence
\[  0\to K\to F \to \bigoplus_{\dim\sigma = m} G_\sigma \to 0.\]
Here $K = F|_{\Delta^{\leq m-1}}$, hence $K$ is semi-Gorenstein on
$\Delta^{\leq m-1}$. The sheaves $G_\sigma$ are clearly
semi-Gorenstein on $\Delta^{\leq m}$. For $\tau \in \Delta^{\leq m-1}$
we get an exact sequence of complexes
\[ 0\to \Cd_m(\Star\tau,K) \to \Cd_m(\Star\tau,F) \to
\bigoplus_{\dim\sigma = m}\Cd_m(\Star\tau,G_\sigma) \to 0,\]
hence an exact sequence of cohomology
\[ 0\to H^0(\Cd_m(\Star\tau,F)) \to \bigoplus_{\dim\sigma = m}
H^0(\Cd_m(\Star\tau,G_\sigma))\to  H^0(\Cd_{m-1}(\Star\tau,K)) \to
0.\]
Thus, we get an exact sequence of sheaves
\[ 0\to K^\vee \to \bigoplus_{\dim\sigma = m} G_\sigma^\vee \to
F^\vee\to 0,\]
where $K^\vee$ is constructed as a sheaf on $\Delta^{\leq m-1}$. By
induction on $m$ we may assume that $K^\vee$ is semi-Gorenstein on
$\Delta^{\leq m-1}$. The sheaves $G_\sigma^\vee$ are constant sheaves
on $[\sigma]$ extended by zero, hence semi-Gorenstein on $\Delta^{\leq
  m}$. Now a long-exact cohomology sequence 
\[ \ldots \to \bigoplus_{\dim\sigma = m}
H^i(\Cd_m(\Star\tau,G_\sigma^\vee)) \to H^i(\Cd_m(\Star\tau,F^\vee))
\to H^i(\Cd_{m-1}(\Star\tau,K^\vee)) \to \ldots \]
shows that $F^\vee$ is also semi-Gorenstein. \qed

\begin{lemma} The functor $(\cdot)^\vee$ is an anti-involution on the
  category of semi-Gorenstein sheaves on $\Delta^{\leq m}$. 
\end{lemma}

{\bf Proof.} We construct a canonical isomorphism $F^{\vee\vee}\isom
F$ on $\Delta^{[m-1,m]}$. If 
$\dim\sigma = m$ then we have canonically
$(F^{\vee\vee})_\sigma = F_\sigma^{**} \isom F_\sigma$. We claim that
this isomorphism extends uniquely to an isomorphism
$F^{\vee\vee}\isom F$. For $\dim\tau=m-1$, $F_\tau^\vee$ is defined by
the exact sequence
\[ 0\to  (F_\tau^\vee)^* \to \bigoplus_{\sigma> \tau} F_\sigma \to
F_\tau \to 0.\]
Similarly, $(F^{\vee\vee})_\tau$ is defined by
\[ 0\to  ((F^{\vee\vee})_\tau)^* \to \bigoplus_{\sigma> \tau} F^\vee_\sigma \to
F^\vee_\tau \to 0.\]
Since the second sequence is the dual of the first, this gives the
isomorphism $(F^{\vee\vee})_\tau \isom F_\tau$ compatible with
restrictions $res^\sigma_\tau$. \qed

Let us now return to sheaves of $\cA$-modules on $\Delta^{\leq m}$. We
will be considering sheaves of the type $\cL\otimes F$, where $F$ is a
semi-Gorenstein sheaf and the tensor product is over $\RR$.

\begin{proposition} \label{prop-dual} Let $F$ be a semi-Gorenstein sheaf on $\Delta^{\leq m}$.
Then in  $D^b(A_m)$ we have an isomorphism 
\[ \Cd_m(\Delta^{\leq m},\cL\otimes F) \isom RHom
  (\Cd_m(\Delta^{\leq m},\cL\otimes F^\vee),\omega_m).\]
\end{proposition}

{\bf Proof.} Bressler and Lunts \cbl construct a dualizing functor $\DD$ on the
derived category $D^b(\cM(\cA))$ of finitely generated $\cA$-modules
and show that 
\[ R\Gamma(\Delta^{\leq m},\DD(\cF))\isom  RHom (\Cd_m(\Delta^{\leq
  m},\cF),\omega_m) \] 
for any $\cF\in D^b(\cM(\cA))$. This result together with
$\DD(\cL)\isom \cL$ gives duality between $\Gamma(\Delta^{\leq
  m},\cL)$ and $\Cd_m(\Delta^{\leq m},\cL)$, hence \Po duality
when $m=n$. 

We show that for $\cF\in D^b(\cM(\cA))$ there exists a canonical
isomorphism 
\begin{equation} \label{eq-dual}
 \Cd_m(\Delta^{\leq  m},\DD(\cF)\otimes F) \isom RHom( \Cd_m(\Delta^{\leq
  m},\cF\otimes F^\vee),\omega_m).
\end{equation} 
Then $\DD(\cL)\isom \cL$ proves the proposition.

The complex of sheaves $\DD(\cF)$ has stalk at $\sigma$ isomorphic to
\[ RHom_{D^b(A_m)}(\Cd_m([\sigma],\cF),\omega_m).\]
To make this an actual complex of sheaves, one first resolves $\cF$ in
projective sheaves of free $A_m$-modules of the type $\cP_\sigma$,
which is the constant sheaf $A_m$ on $[\sigma]$, extended by zero to a
sheaf on $\Delta^{\leq m}$; then $\Cd_m([\sigma],\cP_\sigma)$ is a
complex of free $A_m$-modules and we can replace $RHom_{D^b(A_m)}$ with
$Hom_{A_m}$. Now it suffices to prove the isomorphism
(\ref{eq-dual}) for $\cF = \cP_\sigma$.

Let us fix a $\sigma\in\Delta^{\leq m}$ and consider the double complex
of free $A_m$-modules 
\[ ( 0\to K^{-m}\to K^{-m+1}\to\ldots\to K^0\to 0)\otimes \omega_m,\]
where
\[ K^{-m+i} = \bigoplus_{\tau\leq \sigma,\; \dim\tau = i} \Cd_m (\Star\tau,
F),\]
and the differentials are induced by the projections with signs $or^\sigma_\tau$. We claim that this double
complex represents both sides of the isomorphism  (\ref{eq-dual}) when
$\cF = \cP_\sigma$. 

For a cone $\rho\in \Delta^{\leq m}$ the stalk of
$\DD(\cP_\sigma)\otimes F$ is the complex
\begin{align*} (\DD(\cP_\sigma)\otimes F)_\rho &=
  Hom_{A_m}(\Cd_m([\rho],\cP_\sigma),\omega_m) \otimes F_\rho \\
&= Hom_{A_m}(\Cd_m([\rho\cap\sigma],A_m),\omega_m) \otimes F_\rho \\
&= (\Cd_m([\rho\cap\sigma],\RR))^*\otimes \omega_m \otimes F_\rho.
\end{align*}
This is precisely the contribution from $\rho$ to the double
complex: $\rho\in\Star\tau$ for $\tau\leq \sigma$ if and only if
$\tau\in [\rho\cap\sigma]$. 

 On the other hand, taking cohomology inside each $K^i$, the
double complex becomes
\begin{align*} ( (F_0^\vee)^* \to \bigoplus_{\tau\leq \sigma, \;
  \dim\tau=1} (F_\tau^\vee)^* \to\ldots) \otimes \omega_m &=
  \Cd_m([\sigma],F^\vee)^* \otimes \omega_m \\
&= Hom_{A_m}(\Cd_m([\sigma],F^\vee\otimes A_m), \omega_m) \\
&= Hom_{A_m}(\Cd_m(\Delta^{\leq m},F^\vee\otimes \cP_\sigma),
  \omega_m)
\end{align*}
which is the right hand side of (\ref{eq-dual}). \qed

 To simplify notation, let us denote 
\[ P_m(F) = P_m( \Cd_m(\Delta^{\leq m},\cL\otimes F)).\]
Then for $F$ a semi-Gorenstein sheaf on $\Delta^{\leq m}$, the
proposition gives that $P_m(F)$ and $P_m(F^\vee)$ are dual to each
other. If $F$ is Gorenstein then $P_m(F) = P_m(F^\vee)$.

 \begin{lemma}\label{lem-G} Let $F$ be a semi-Gorenstein sheaf on $\Delta^{\leq
    m}$. Assume that
\begin{enumerate}
\item We have a surjective morphism $F^\vee\to F$ with kernel $G$.
\item $P_m(F)$ is a $cd$-polynomial 
\[   P_m(F) = f(c,d) +g(c,d),\]
where $f$ and $g$ are homogeneous $cd$-polynomials of degree $m$ and
$m-1$, respectively.
\end{enumerate}
Then $G$ is a Gorenstein sheaf on $\Delta^{\leq m-1}$ and 
\[ P_{m-1}(G) = g(c,d).\]
\end{lemma}

{\bf Proof.}
From the exact sequence
\[ 0 \to G \to F^\vee \to F\to 0\]
we have $\dim G_\sigma = \dim F^\vee_\sigma-\dim F_\sigma$. Computing
the cohomology long-exact sequence of 
\[ 0\to \Cd_m(\Star\sigma, G) \to \Cd_m(\Star\sigma, F^\vee) \to
\Cd_m(\Star\sigma, F) \to 0,\]
we get that $H^i( \Cd_m(\Star\sigma, G)) = 0$ for $i>1$, hence $G$ is
semi-Gorenstein on $\Delta^{\leq m-1}$, and 
\[ 0\to H^0(\Cd_m(\Star\sigma, F^\vee)) \to H^0(\Cd_m(\Star\sigma, F))
\to  H^1( \Cd_m(\Star\sigma, G)) \to 0\]
shows that 
\[ \dim G^\vee_\sigma = \dim F^\vee_\sigma - \dim F_\sigma = \dim
G_\sigma,\]
thus $G$ is Gorenstein.

Proposition~\ref{prop-dual} gives that $P_m(F^\vee) =
f(c,d)+g(c,d)t_m$. Now from the exact sequence 
\[ 0 \to G\otimes \cL \to F^\vee \otimes\cL\to F\otimes\cL\to 0,\]
we get 
\begin{align*} P_{m-1}(G)(1-t_n) = P_m(G) &= -P_m(\Cd_m(\Delta^{\leq
    m},\cL\otimes G)) \\  &= P_m(F)-P_m(F^\vee) = g(c,d)(1-t_n).
\end{align*} 
This proves the last statement. \qed

To apply the previous lemma, we need to see when a surjective morphism
$F^\vee\to F$ exists. Since $F^\vee$ is semi-Gorenstein on $\Delta^{\leq
    m}$, it suffices to construct such a morphism on
$\Delta^{[m-1,m]}$. Moreover, such a morphism is uniquely determined
if we know it on $\Delta^{m}$. It is also easy to see that such a
morphism is surjective if and only if its restriction to $\Delta^m$ is
surjective. Let us see when a morphism defined on  $\Delta^{m}$
extend to a morphism on  $\Delta^{[m-1,m]}$. 

Suppose $G$ and $F$ are semi-Gorenstein sheaves on $\Delta^{\leq m}$,
and let $\phi_\sigma: G_\sigma\to F_\sigma$ be given for all
$\sigma\in\Delta^m$. This map extends to a morphism of sheaves $G\to
F$ if for all $\tau\in\Delta^{m-1}$ the map $\oplus \phi_\sigma$ takes
\[ \Ker( \oplus_{\sigma>\tau} G_\sigma
\stackrel{res^\sigma_\tau}{\longrightarrow} G_\tau) \to
\Ker(\oplus_{\sigma>\tau} F_\sigma 
\stackrel{res^\sigma_\tau}{\longrightarrow} F_\tau).\]
Indeed, then we can lift an element in $G_\tau$ to
$\oplus_{\sigma>\tau} G_\sigma$ and map it to $\oplus_{\sigma>\tau}
F_\sigma$ followed by restriction to $F_\tau$.

When $G=F^\vee$ then $G_\sigma = F_\sigma^*$ and 
\[ G_\tau =  \Ker(\oplus_{\sigma>\tau} F_\sigma
\stackrel{or^\sigma_\tau res^\sigma_\tau}{\longrightarrow} F_\tau)^* =
\Coker (F_\tau^*  \stackrel{or^\sigma_\tau
  (res^\sigma_\tau)^*}{\longrightarrow} \oplus_{\sigma>\tau} F_\sigma^*).\]
Now given $\phi_\sigma: F_\sigma^* \to F_\sigma$ for
all $\sigma\in\Delta^m$, the compatibility condition becomes that the
composition 
\begin{equation}\label{eq-comp} F_\tau^* \stackrel{or^\sigma_\tau
  (res^\sigma_\tau)^*}{\longrightarrow} \oplus_{\sigma>\tau} F_\sigma^*
  \stackrel{\phi_\sigma}{\longrightarrow} \oplus_{\sigma>\tau} F_\sigma 
\stackrel{res^\sigma_\tau}{\longrightarrow} F_\tau
\end{equation}
must be zero.

As an example, take $F$ to be the constant sheaf $\RR$ on
$\Delta^{\leq m}$. Then $\phi_\sigma: \RR^*\to \RR$ is given by a constant
$c_\sigma$ and the sequence (\ref{eq-comp}) becomes
\[ \RR_\tau^* \stackrel{or^\sigma_\tau}{\longrightarrow} \oplus_{\sigma>\tau} \RR_\sigma^*
  \stackrel{c_\sigma}{\longrightarrow} \oplus_{\sigma>\tau} \RR_\sigma 
\stackrel{1}{\longrightarrow} \RR_\tau,\]
the composition of which is 
\[ \sum_{\sigma>\tau} or^\sigma_\tau c_\sigma.\]
This sum is zero for all $\tau$ precisely when the collection
$(c_\sigma)$ defines a class in $H^0(\Cd(\Delta^{\leq m},\RR))$. When
$m<n$ then such a class is defined as follows. Choose $d_\pi\in\RR$
for each $\pi\in\Delta^{m+1}$, and set $c_\sigma = \sum_{\pi>\sigma}
or^\pi_\sigma d_\pi$. Then for $\tau\in\Delta^{m-1}$ we have
\[  \sum_{\sigma>\tau} or^\sigma_\tau c_\sigma =
\sum_{\pi>\sigma>\tau} or^\sigma_\tau or^\pi_\sigma d_\pi = \sum_\pi
d_\pi ( \sum_{\pi> \sigma>\tau}  or^\sigma_\tau or^\pi_\sigma) = 0.\]
Choosing $d_\pi$ generically ensures that all $c_\sigma$ are nonzero,
hence $\phi: F^\vee\to F$ is surjective.

\begin{proposition} \label{prop-surj}
Let $F$ be a semi-Gorenstein sheaf on $\Delta^{\leq m+1}$. Then there
exists a surjective morphism 
\[ F|_{\Delta^{\leq m}}^\vee \to F|_{\Delta^{\leq m}}.\]
\end{proposition}

{\bf Proof.} We look for $\phi_\sigma: F_\sigma^* \to F_\sigma$ for
$\sigma\in\Delta^m$ in the form of a composition
\[ \phi_\sigma: F_\sigma^* \stackrel{or^\pi_\sigma
  (res^\pi_\sigma)^*}{\longrightarrow} \oplus_{\pi>\sigma} F_\pi^*
  \stackrel{\phi_\pi}{\longrightarrow} \oplus_{\pi>\sigma} F_\pi 
\stackrel{res^\pi_\sigma}{\longrightarrow} F_\sigma,\]
where $\phi_\pi: F_\pi^*\to F_\pi$ are defined for all
  $\pi\in\Delta^{m+1}$ independently of $\sigma$. Plugging this
  sequence into the sequence (\ref{eq-comp}), we get for each $\tau$
\[ F_\tau^* \stackrel{or^\sigma_\tau
  (res^\sigma_\tau)^*}{\longrightarrow} \oplus_{\sigma>\tau}
F_\sigma^* \stackrel{or^\pi_\sigma
  (res^\pi_\sigma)^*}{\longrightarrow} \oplus_{\pi>\sigma>\tau} F_\pi^*
  \stackrel{\phi_\pi}{\longrightarrow} \oplus_{\pi>\sigma>\tau} F_\pi 
\stackrel{res^\pi_\sigma}{\longrightarrow}  \oplus_{\sigma>\tau} F_\sigma 
\stackrel{res^\sigma_\tau}{\longrightarrow} F_\tau,\]
whose composition is zero because for a fixed $\pi$ and $\tau$ we have
\[\sum_{\pi>\sigma>\tau} or^\sigma_\tau or^\pi_\sigma res^\sigma_\tau
res^\pi_\sigma = 0.\] 

Thus, it suffices to find for all  $\pi\in\Delta^{m+1}$ a linear map
$\phi_\pi: F_\pi^*\to F_\pi$ such that the composition $\phi_\sigma$
is an isomorphism for each $\sigma\in\Delta^m$. We claim that a
general collection of $\phi_\pi$ satisfies this condition: for each
$\sigma$ there is a proper algebraic subset in $\oplus_\pi
\Hom_\RR( F_\pi^*,F_\pi)$ for which $\phi_\sigma$ is not an
isomorphism. 

Let us fix one $\sigma$. The condition for $\phi_\sigma$ not being an
isomorphism is clearly 
algebraic. Thus, it suffices to find one collection $\phi_\pi$ that
gives an isomorphism $\phi_\sigma$. We may assume that
$or^\pi_\sigma=1$, otherwise replace $\phi_\pi$ by $-\phi_\pi$. Now
an isomorphism $\phi_\pi$ is given by a non-degenerate bilinear form
$B_\pi$ on $F_\pi^*$ and $\phi_\sigma$ is given by the restriction of
$\oplus_\pi B_\pi$ to the subspace $F_\sigma^* \subset \oplus_\pi
F_\pi^*$, provided that this restriction is non-degenerate. If we take
$B_\pi$ to be positive definite symmetric bilinear forms, then their
restriction is 
also positive definite, hence defines an isomorphism $\phi_\sigma$. \qed

\begin{lemma} Let $F$ be a semi-Gorenstein sheaf on $\Delta^{\leq m+1}$ such
  that $P_{m+1}(F)= f(c,d)+g(c,d)$, where $f$ and $g$ are homogeneous
  $cd$-polynomials of degree $m+1$ and $m$, respectively. If
\[ f(c,d) = f_m(c,d)c+g_m(c,d)d,\]
then 
\[ P_{m}(F|_{\Delta^{\leq m}}) = f_m(c,d)+g_m(c,d).\]
\end{lemma}

{\bf Proof.} Consider the exact sequence
\[ 0\to \Cd_m(\Delta^{\leq m}, F\otimes \cL)[-1] \to
\Cd_{m+1}(\Delta^{\leq m+1}, F \otimes \cL) \to
\bigoplus_{\dim\sigma = m+1} F_\sigma\otimes \cL_\sigma \to 0.\]
Let $P_{m+1}(\oplus_{\dim\sigma = m+1} F_\sigma\otimes \cL_\sigma) =
h(c,d)$, which is a homogeneous $cd$-polynomial of degree $m$. Then
\[ -P_m(F|_{\Delta^{\leq m}})(1-t_{m+1}) - P_{m+1}(F) + h(c,d) = 0.\]
It follows that $P_m(F|_{\Delta^{\leq m}})$ is the
coefficient of $t_{m+1}$ in  $P_{m+1}(F)$. \qed

Let us define operations $\cC$ and $\cD$ on semi-Gorenstein
sheaves. If $F$ is a semi-Gorenstein sheaf on $\Delta^{\leq m}$, let $\cC(F)$
be its restriction to $\Delta^{\leq m-1}$, and let $\cD(F)$ be the
sheaf $G$ on $\Delta^{\leq m-2}$ from Lemma~\ref{lem-G} and
Proposition~\ref{prop-surj}. Such a sheaf $G$
is not uniquely defined, however, the dimensions of its stalks are.
Given a $cd$-monomial $w(c,d)$ of degree $m$, we can apply the
operation $w(\cC,\cD)$ to a semi-Gorenstein sheaf $F$ on $\Delta^{\leq
  m}$. The result is a semi-Gorenstein sheaf on the cone $0$, hence a
vector space. 

\begin{theorem}\label{thm-final} Let $F$ be a semi-Gorenstein sheaf on
  $\Delta^{\leq 
    m}$ such that $P_m(F) = f(c,d) +g(c,d)$, where $f$ and $g$ are
    homogeneous $cd$-polynomials of degree $m$ and $m-1$, respectively. 
    If $w(c,d)$ is a $cd$-monomial of degree $m$, then 
\[ \dim w(\cC,\cD)(F)_0 \]
is the coefficient of $w(c,d)$ in $f(c,d)$.
\end{theorem}

{\bf Proof.}
We prove it by induction on $m$, the case $m=0$ being trivial. 

Let $f = f_{m-1}(c,d)c+ g_{m-1}(c,d)d$. Then by the previous lemma
\[ P_{m-1}(\cC(F)) =  f_{m-1}(c,d)+ g_{m-1}(c,d),\]
and by Lemma~\ref{lem-G}
\[ P_{m-2}(\cD(F)) = g_{m-1}(c,d).\]
Now we are done by induction: If $w=w_{m-1}c$ then the coefficient of
$w$ in $f(c,d)$ is the coefficient of $w_{m-1}$ in $f_{m-1}(c,d)$,
which is $\dim w_{m-1}(\cC,\cD)(\cC(F))_0$. If  $w=w_{m-2}d$ then the
coefficient of $w$ in $f(c,d)$ is the coefficient of $w_{m-2}$ in
$g_{m-1}(c,d)$, which is $\dim w_{m-2}(\cC,\cD)(\cD(F))_0$. 
\qed

Applying this to $F=\RR_\Delta$, we get

\begin{corollary} If $\Delta$ is a complete $n$-dimensional fan, then
  $P_n(\Delta)$ is a 
  homogeneous $cd$-polynomial with non-negative integer coefficients. \qed
\end{corollary}

\begin{corollary} If $\Delta$ is a quasi-convex $n$-dimensional fan,
  then $P_n(\Delta)$ is a  $cd$-polynomial with non-negative integer
  coefficients. 
\end{corollary}

{\bf Proof.} Embed $\Delta$ in a complete fan $\overline{\Delta}$ of
 the same dimension.  Let $F$ be
the constant sheaf $\RR$ on $\Delta\setmin\partial\Delta$, extended by
zero to a sheaf on $\overline{\Delta}$. Then
\[ P_n(F) = P_n(\Delta) = f(c,d)+P_{n-1}(\partial\Delta).\]
Now $F$ is semi-Gorenstein, in fact $F = \RR_\Delta^\vee$, where
$\RR_\Delta$ is the constant sheaf $\RR$ on $\Delta$ extended by zero. 
The theorem gives that $f(c,d)$ has non-negative integer coefficients,
the previous corollary proves the same for $P_{n-1}(\partial\Delta)$.\qed

\subsection{Gorenstein* posets} 
Let $\Lambda$ be a Gorenstein* poset. Then $B(\Lambda)$ is a
well-complemented Gorenstein lattice and the theory of semi-Gorenstein
sheaves on $B(\Lambda)$ 
works the same way as for complete fans. When working with
$\cB$-modules on $B(\Lambda)$, the difference is that restriction to
the $m$-skeleton should be restriction to $B(\Lambda^{\leq m})$,
rather than restriction to $B(\Lambda)^{\leq m}$. Then for
Lemma~\ref{lem-G} to hold, we need the sheaf $G$ to be supported on
$B(\Lambda^{\leq m})$. 

If $F$ is a sheaf of finite dimensional vector spaces on
$B(\Lambda^{\leq m})$, let $\dim F:B(\Lambda^{\leq m})\to\ZZ$ be
the element in $I(B(\Lambda^{\leq m}))$ such that 
\[ \dim F (\sigma) = \dim F_\sigma.\]
We consider sheaves $F$ on $B(\Lambda^{\leq m})$ satisfying the
condition 
\begin{equation} \label{eq-pull}
 \dim F = \pi^*(f) 
\end{equation}
for some $f\in I(\Lambda^{\leq m})$.
Recall that $\pi: B(\Lambda) \to \Lambda$ is the projection.

\begin{lemma} If $F$ is a semi-Gorenstein sheaf on $B(\Lambda^{\leq
    m})$ satisfying condition~(\ref{eq-pull}), then $F^\vee$ satisfies
    the same condition.
\end{lemma}

{\bf Proof.} Recall the definition of the operation $E$ given in the
introduction. We have 
\[ \dim F^\vee = E(\dim F).\]
Now it suffices to prove that $E$ commutes with $\pi^*$:
\[ E\circ\pi^* (f) = \pi^* \circ E (f).\]
Then $\dim F^\vee = \pi^*(E(f))$. 

For $\sigma\in  B(\Lambda^{\leq m})$, 
\begin{align*} E\circ\pi^* (f)(x) 
&= \sum_{y\geq x} (-1)^{m-\dim y} f(\pi(y))\\ 
&= \sum_{\sigma\in\Lambda,\; \sigma\geq \pi(x)} (-1)^{m-\dim \sigma}
  f(\sigma)   
\sum_{y\geq x,\;
  \pi(y) = \sigma} (-1)^{\dim \sigma-\dim y}.
\end{align*}
The last sum runs over all elements in $[\sigma,{\bf 1})\subset
  B([{\bf 0},x])$. Since $B([{\bf 0},x])$ is Eulerian, this sum is $1$. \qed

The previous lemma implies that if $F$ satisfies
condition~(\ref{eq-pull}), then the sheaf $G$ constructed in
Lemma~\ref{lem-G} again satisfies the same condition. It follows that 
the operations $\cC$ and $\cD$ are well-defined on such sheaves and
the results also satisfy the same condition. Now all proofs above
work for semi-Gorenstein sheaves $F$ on $B(\Lambda^{\leq m})$
satisfying condition~(\ref{eq-pull})  with $\cL$ replaces by $\cB$.

This proves Theorem~\ref{thm-main}. The operations $\cC$ and $\cD$ on
semi-Gorenstein sheaves correspond to operations $C$ and $D$ defined
in the introduction. Thus, 
Theorem~\ref{thm-final} implies Proposition~\ref{prop-constr} for not
only complete fans but also for Gorenstein* posets.

\end{document}